\documentclass[12pt]{article}

\usepackage{amsmath, epsfig, cite}
\usepackage{amssymb}
\usepackage{amsfonts}
\usepackage{latexsym}
\usepackage{float}
\usepackage{color}
\usepackage{booktabs}
\usepackage{graphicx}
\usepackage{caption}

\newtheorem{thm}{Theorem}[section]

\newtheorem{lem}[thm]{Lemma}

\numberwithin{equation}{section}

\newcommand{\qed}{{\hfill$\square$}\medskip}
\UseRawInputEncoding
\begin{document}

\begin{center}
{\large\bf An extension of Gauss congruences for Ap\'ery numbers}
\end{center}

\vskip 2mm \centerline{Ji-Cai Liu}
\begin{center}
{\footnotesize Department of Mathematics, Wenzhou University, Wenzhou 325035, PR China\\[5pt]
{\tt jcliu2016@gmail.com} \\[10pt]
}
\end{center}

\vskip 0.7cm \noindent{\bf Abstract.}
Osburn, Sahu and Straub introduced the numbers:
\begin{align*}
A_n^{(r,s,t)}=\sum_{k=0}^n{n\choose k}^r{n+k\choose k}^s{2k\choose n}^t,
\end{align*}
for non-negative integers $n,r,s,t$ with $r\ge 2$, which includes two kinds of Ap\'ery numbers and four kinds of Ap\'ery-like numbers as special cases, and showed that the numbers $\{A_n^{(r,s,t)}\}_{n\ge 0}$ satisfy the Gauss congruences of order $3$. We establish an extension of Osburn--Sahu--Straub congruence through Bernoulli numbers, which is one step deep congruence of the Gauss congruence for $A_n^{(r,s,t)}$.


\vskip 2mm
\noindent{\it MR Subject Classifications}: 11B50, 11B65, 11B68

\section{Introduction}
In 1979, Ap\'ery \cite{apery-asterisque-1979} introduced the two sequences $\{a_n\}_{n\ge 0}$ and $\{b_n\}_{n\ge 0}$ through $3$-term recurrences in the proof of the irrationality of $\zeta(3)$ and $\zeta(2)$:
\begin{align}
&(n+1)^3 a_{n+1}- (2n+1)(17n^2+17n+5)a_n+n^3 a_{n-1}=0,\quad (a_0=1,a_1=5),\label{a-1}\\[5pt]
&(n+1)^2b_{n+1}-(11n^2+11n+3)b_n-n^2b_{n-1}=0,\quad(b_0=1,b_1=3).\label{a-2}
\end{align}
The two sequences $\{a_n\}_{n\ge 0}$ and $\{b_n\}_{n\ge 0}$ are known as the famous Ap\'ery numbers, which possess the binomial sum formulae:
\begin{align*}
&a_n=\sum_{k=0}^n{n\choose k}^2{n+k\choose k}^2,\\[5pt]
&b_n=\sum_{k=0}^n{n\choose k}^2{n+k\choose k}.
\end{align*}

Zagier \cite{zagier-b-2009} considered the following recurrence related to \eqref{a-2}:
\begin{align}
(n+1)^2u_{n+1}-(An^2+An+\lambda)u_n+Bn^2u_{n-1}=0,\quad (u_{-1}=0,u_0=1),\label{rec-apery-1}
\end{align}
and searched for triples $(A,B,\lambda)\in \mathbb{Z}^3$ such that the solution of the recurrence \eqref{rec-apery-1} is an integer sequence $\{u_n\}_{n\ge 0}$. Six sporadic sequences are found in Zagier's search, which include the desired solution $\{b_n\}_{n\ge 0}$. The six sporadic sequences are listed in the following table.
\begin{table}[H]
\centering
\scalebox{0.8}{
\begin{tabular}{cccc}
\toprule
$(A,B,\lambda)$&Name&Other names&Formula\\
\midrule
$(7,-8,2)$&{\bf A}&Franel numbers&$u_n=\sum_{k=0}^n{n\choose k}^3$ \\[7pt]
$(9,27,3)$&{\bf B}&&$u_n=\sum_{k=0}^n(-1)^k3^{n-3k}{n\choose 3k}{3k\choose 2k}{2k\choose k}$ \\[7pt]
$(10,9,3)$&{\bf C}&&$u_n=\sum_{k=0}^n{n\choose k}^2{2k\choose k}$\\[7pt]
$(11,-1,3)$&{\bf D}&Ap\'ery numbers&$u_n=\sum_{k=0}^n{n\choose k}^2{n+k\choose k}$\\[7pt]
$(12,32,4)$&{\bf E}&&$u_n=\sum_{k=0}^n4^{n-2k}{n\choose 2k}{2k\choose k}^2$\\[7pt]
$(17,72,6)$&{\bf F}&&$u_n=\sum_{k=0}^n(-1)^k8^{n-k}{n\choose k}\sum_{j=0}^k{k\choose j}^3$\\
\bottomrule
\end{tabular}
}
\end{table}

Almkvist and Zudilin \cite{az-b-2006} studied the other recurrence related to \eqref{a-1}:
\begin{align}
(n+1)^3u_{n+1}-(2n+1)(an^2+an+b)u_n+cn^3u_{n-1}=0,\quad (u_{-1}=0,u_0=1),\label{a-3}
\end{align}
and searched for triples $(a,b,c)\in \mathbb{Z}^3$ such that the solution of the recurrence \eqref{a-3}
is an integer sequence $\{u_n\}_{n\ge 0}$. They also found six sporadic sequences, which include the desired solution $\{a_n\}_{n\ge 0}$.

Cooper \cite{cooper-rj-2012} considered a more general recurrence:
\begin{align}
(n+1)^3u_{n+1}-(2n+1)(an^2+an+b)u_n+n(cn^2+d)u_{n-1}=0,\label{rec-apery-2}
\end{align}
with $u_{-1}=0$ and $u_0=1$. Note that the case $d=0$ of \eqref{rec-apery-2} reduces to \eqref{a-3}. He found three additional sporadic sequences $s_7,s_{10}$ and $s_{18}$. The nine sporadic sequences are listed in the following table.

\begin{table}[H]
\centering
\scalebox{0.8}{
\begin{tabular}{cccc}
\toprule
$(a,b,c,d)$&Name&Other names&Formula\\
\midrule
$(7,3,81,0)$&$(\delta)$&Almkvist-Zudilin numbers&$u_n=\sum_{k=0}^n(-1)^k3^{n-3k}{n\choose 3k}{n+k\choose k}{3k\choose 2k}{2k\choose k}$ \\[7pt]
$(11,5,125,0)$&$(\eta)$&&$u_n=\sum_{k=0}^n(-1)^k{n\choose k}^3({4n-5k-1\choose 3n}+{4n-5k\choose 3n})$ \\[7pt]
$(10,4,64,0)$&$(\alpha)$&Domb numbers&$u_n=\sum_{k=0}^n{n\choose k}^2{2k\choose k}{2n-2k\choose n-k}$ \\[7pt]
$(12,4,16,0)$&$(\epsilon)$&&$u_n=\sum_{k=0}^n{n\choose k}^2{2k\choose n}^2$\\[7pt]
$(9,3,-27,0)$&$(\zeta)$&&$u_n=\sum_{k=0}^n\sum_{l=0}^n{n\choose k}^2{n\choose l}{k\choose l}{k+l\choose n}$ \\[7pt]
$(17,5,1,0)$&$(\gamma)$&Ap\'ery numbers&$u_n=\sum_{k=0}^n{n\choose k}^2{n+k\choose k}^2$ \\[7pt]
$(13,4,-27,3)$&$s_7$&&$u_n=\sum_{k=0}^n{n\choose k}^2{n+k\choose k}{2k\choose n}$ \\[7pt]
$(6,2,-64,4)$&$s_{10}$&Yang-Zudilin numbers&$u_n=\sum_{k=0}^n{n\choose k}^4$ \\[7pt]
$(14,6,192,-12)$&$s_{18}$&&$u_n=\sum_{k=0}^n(-1)^k{n\choose k}{2k\choose k}{2n-2k\choose n-k}({2n-3k-1\choose n}+{2n-3k\choose n})$ \\
\bottomrule
\end{tabular}
}
\end{table}

An integer sequence $\{u_n\}_{n\ge 1}$ is said to satisfy the Gauss congruences of order $r$ if
$u_{np^k}\equiv u_{np^{k-1}} \pmod{p^{rk}}$ for all positive integers $n,k$ and all primes $p\ge r+1$.
Let $\mathbb{N}$ denote the set of non-negative integers, $\mathbb{Z}^{+}$ denote the set of positive integers, and $\mathbb{Q}$ denote the set of rational numbers.

As early as 1982, Gessel \cite{gessel-jnt-1982} proved that for all primes $p\ge 5$ and $n\in \mathbb{Z}^{+}$,
\begin{align*}
a_{np}\equiv a_{n} \pmod{p^3}.
\end{align*}
Coster \cite{coster-phd-1988} further showed that the Ap\'ery numbers $\{a_n\}_{n\ge 0}$ satisfy the Gauss congruences of order $3$, namely,
\begin{align*}
a_{np^m}\equiv a_{np^{m-1}}\pmod{p^{3m}},
\end{align*}
for primes $p\ge 5$ and $n,m\in \mathbb{Z}^{+}$.

In 2016, Osburn, Sahu and Straub \cite{oss-pems-2016} defined
\begin{align*}
A_n^{(r,s,t)}=\sum_{k=0}^n{n\choose k}^r{n+k\choose k}^s{2k\choose n}^t,
\end{align*}
for $n,r,s,t\in \mathbb{N}$ with $r\ge 2$. Note that $A_n^{(r,s,t)}$ includes the sequences ${\bf A}, {\bf D}, (\epsilon),(\gamma), s_7, s_{10}$ as special cases.
\begin{table}[H]
\caption*{\text{Special cases of $A_n^{(r,s,t)}$}}
\centering
\scalebox{0.8}{
\begin{tabular}{cccc}
\toprule
$(r,s,t)$&Name&Other names&Formula\\
\midrule
$(3,0,0)$& {\bf A} & Franel numbers & $u_n=\sum_{k=0}^n{n\choose k}^3$\\[7pt]
$(2,1,0)$& {\bf D}&Ap\'ery numbers & $u_n=\sum_{k=0}^n{n\choose k}^2{n+k\choose k}$\\[7pt]
$(2,0,2)$&$(\epsilon)$& &$u_n=\sum_{k=0}^n{n\choose k}^2{2k\choose n}^2$\\[7pt]
$(2,2,0)$&$(\gamma)$&Ap\'ery numbers&$u_n=\sum_{k=0}^n{n\choose k}^2{n+k\choose k}^2$ \\[7pt]
$(2,1,1)$&$s_7$&&$u_n=\sum_{k=0}^n{n\choose k}^2{n+k\choose k}{2k\choose n}$\\[7pt]
$(4,0,0)$&$s_{10}$&Yang-Zudilin numbers&$u_n=\sum_{k=0}^n{n\choose k}^4$\\
\bottomrule
\end{tabular}
}
\end{table}

Osburn, Sahu and Straub \cite{oss-pems-2016} also showed that the numbers $\{A_n^{(r,s,t)}\}_{n\ge 0}$ also satisfy the Gauss congruences of order $3$, namely,
\begin{align}
A_{np^m}^{(r,s,t)}\equiv A_{np^{m-1}}^{(r,s,t)} \pmod{p^{3m}},\label{oos-cong}
\end{align}
for all primes $p\ge 5$, $n,m\in \mathbb{Z}^{+}$ and $r,s,t\in \mathbb{N}$ with $r\ge 2$.

The Bernoulli numbers $\{B_n\}_{n\ge 0}$ are defined by the generating function:
\begin{align*}
\frac{z}{e^z-1}=\sum_{n=0}^{\infty}B_n \frac{z^n}{n!}.
\end{align*}
In 2020, Sun \cite[Conjectures 5.1 and 5.3]{sunzh-a-2020} conjectured a series of supercongruences
concerning $u_p,u_{2p}$ and $u_{3p}$ modulo $p^4$, where $\{u_n\}_{n\ge 0}$ is one of the sequences ${\bf A},{\bf D},(\alpha),(\gamma),(\epsilon)$.

The motivation of the paper is to establish an extension of the Gauss congruence
\eqref{oos-cong}, which involves the Bernoulli numbers. The main result is stated as follows.

\begin{thm}\label{t-1}
Let $p\ge 5$ be a prime, $n,m\in \mathbb{Z}^{+}$, and $r,s,t\in \mathbb{N}$ with $r\ge 2$.
Then
\begin{align*}
A_{np^m}^{(r,s,t)}\equiv A_{np^{m-1}}^{(r,s,t)}+p^{3m}B_{p-3} \mathcal{A}^{(r,s,t)}_n\pmod{p^{3m+1}},
\end{align*}
where $\mathcal{A}^{(r,s,t)}_n$, independent of $m$ and $p$, are given by
\begin{align*}
\mathcal{A}^{(2,s,t)}_n
&=-\frac{1}{3}\sum_{k=0}^{n} {n\choose k}^2 {n+k\choose k}^s
{2k\choose n}^t nk(sn+sk+2n-2k+4tk-2tn)\\[5pt]
&+\frac{1}{6}\sum_{k=0}^{n-1} {n\choose k}^2{n+k\choose k}^s {2k\choose n}^t (n-k)^2\left(9nt-3ns+24k-18n+14\right)\\[5pt]
&+\frac{1}{6}\sum_{k=0}^{n-1}{n\choose k}^2 {n+k\choose k}^s {2k+1\choose n}^t
(n-k)^2\left(9ns+15nt-24k+6n-10\right),
\end{align*}
\begin{align*}
\mathcal{A}_{n}^{(3,s,t)}
&=-\frac{1}{3}\sum_{k=0}^{n} {n\choose k}^3 {n+k\choose k}^s
{2k\choose n}^t nk(sn+sk+3n-3k+4tk-2tn)\\[5pt]
&+\frac{1}{4}\sum_{k=0}^{n-1}(n-k)^3{n\choose k}^3 {n+k\choose k}^s\left({2k\choose n}^t+{2k+1\choose n}^t\right),
\end{align*}
and
\begin{align*}
&\mathcal{A}_{n}^{(r,s,t)}\\[5pt]
&=-\frac{1}{3}\sum_{k=0}^{n} {n\choose k}^r {n+k\choose k}^s
{2k\choose n}^t nk(sn+sk+rn-rk+4tk-2tn)\quad\text{for $r\ge 4$.}
\end{align*}
\end{thm}

The rest of the paper is organized as follows. We establish some preliminary results in the next section.
The proof of Theorem \ref{t-1} is presented in Section 3, which is divided into two parts.

\section{Preliminaries}
Throughout the paper, let $\sum{'}$ denote the sum over indices not divisible by $p$, $\lfloor x\rfloor$ denote the integral part of real $x$, and $\{k/p^{m}\}$ denote the remainder of $k$ divided by $p^m$, namely,
\begin{align*}
\{k/p^{m}\}=k-p^m\lfloor k/p^{m}\rfloor.
\end{align*}

\begin{lem}(See \cite{granville-b-1997}.)
Let $p\ge 5$ be a prime, and $n,k\in \mathbb{Z}^{+}$ with $n\ge k$. Then
\begin{align}
{np\choose kp}{\Big/}{n\choose k}\equiv 1-\frac{nk(n-k)}{3}p^3B_{p-3} \pmod{p^{\text{ord}_p\left(nk(n-k)\right)+4}},\label{b-1}
\end{align}
where $\text{ord}_p(x)$ denotes the exponent of $p$ in $x\in \mathbb{Q}$.
\end{lem}

\begin{lem}
Let $p\ge 5$ be a prime, $a,b,c,r,s,t\in \mathbb{N}$ and $n,m\in \mathbb{Z}^{+}$. Then
\begin{align}
&{np^m-1\choose a}^r {np^m+b\choose b}^s {c\choose np^m}^t\notag\\[5pt]
&\equiv (-1)^{r(a+\lfloor a/p\rfloor)}{np^{m-1}-1\choose \lfloor a/p\rfloor}^r
{np^{m-1}+\lfloor b/p\rfloor\choose \lfloor b/p\rfloor}^s
{\lfloor c/p\rfloor \choose np^{m-1}}^t\notag\\[5pt]
&\times\left(1-rnp^m\sum_{j=1}^{a}{'}\frac{1}{j}
+snp^m\sum_{j=1}^b{'}\frac{1}{j}+tnp^m\sum_{j=1}^{c}{'}\frac{1}{j}\right)\pmod{p^{m+1}}.\label{b-2}
\end{align}
\end{lem}
{\noindent\it Proof.}
Note that
\begin{align}
{np^m-1\choose a}
&=\prod_{j=1}^a\frac{np^m-j}{j}\notag\\
&=\prod_{\substack{j=1\\ p\nmid j}}^a\frac{np^m-j}{j}\prod_{i=1}^{\lfloor a/p\rfloor}\frac{np^{m-1}-i}{i}\notag\\
&={np^{m-1}-1\choose \lfloor a/p\rfloor}\prod_{\substack{j=1\\ p\nmid j}}^a\frac{np^m-j}{j}\notag\\
&\equiv (-1)^{a+\lfloor a/p\rfloor}{np^{m-1}-1\choose \lfloor a/p\rfloor}
\left(1-np^m\sum_{j=1}^a{'}\frac{1}{j}\right)\pmod{p^{m+1}}.\label{b-3}
\end{align}
By \eqref{b-3}, we have
\begin{align}
{np^m+b\choose b}&=(-1)^b{-np^m-1\choose b}\notag\\[5pt]
&\equiv {np^{m-1}+\lfloor b/p\rfloor\choose \lfloor b/p\rfloor}
\left(1+np^m\sum_{j=1}^b{'}\frac{1}{j}\right)\pmod{p^{m+1}},\label{b-4}
\end{align}
and
\begin{align}
{c\choose np^m}&={np^m+c-np^m\choose c-np^m}\notag\\[5pt]
&\equiv {np^{m-1}+\lfloor c/p\rfloor -np^{m-1}\choose \lfloor c/p\rfloor-np^{m-1}}
\left(1+np^m\sum_{j=1}^{c-np^m}{'}\frac{1}{j}\right)\pmod{p^{m+1}}.\label{b-5}
\end{align}
Since $\sum_{j=1}^{p-1}\frac{1}{j^d} \equiv 0\pmod{p}$ for $d\in \mathbb{Z}^{+}$, we have
\begin{align}
\sum_{j=1}^{a}{'}\frac{1}{j^d}\equiv \sum_{j=1}^{b}{'}\frac{1}{j^d} \pmod{p},
\quad\text{for $d\in \mathbb{Z}^{+}$ and $a\equiv b\pmod{p}$.}\label{b-new-1}
\end{align}
Then we rewrite \eqref{b-5} as
\begin{align}
{c\choose np^m}\equiv {\lfloor c/p\rfloor \choose np^{m-1}}\left(1+np^m\sum_{j=1}^{c}{'}\frac{1}{j}\right) \pmod{p^{m+1}}.\label{b-6}
\end{align}
Combining \eqref{b-3}, \eqref{b-4} and \eqref{b-6}, we arrive at the desired result \eqref{b-2}.
\qed

\begin{lem}
Let $p\ge 5$ be a prime, $n\in \mathbb{N}$ and $m\in \mathbb{Z}^{+}$. Then
\begin{align}
&\sum_{\substack{\lfloor k/p^{m}\rfloor=n\\[3pt] \{k/p^{m}\}<p^{m}/2}}{'}\frac{1}{k^2}\equiv \frac{12n+7}{3}p^m B_{p-3}\pmod{p^{m+1}},\label{b-7}\\
&\sum_{\substack{\lfloor k/p^{m}\rfloor=n\\[3pt] \{k/p^{m}\}>p^{m}/2}}{'}\frac{1}{k^2}\equiv -\frac{12n+5}{3}p^m B_{p-3}\pmod{p^{m+1}}.\label{b-8}
\end{align}
\end{lem}
{\noindent\it Proof.}
By \cite[Theorem 1.2]{sla-camb-2002}, we have
\begin{align}
\sum_{k=1}^{p^m-1}{'}\frac{1}{k^2}\equiv \frac{2}{3}p^mB_{p-3} \pmod{p^{m+1}}.\label{b-9}
\end{align}

Since
\begin{align*}
\frac{1}{(p^m-k)^2}\equiv \frac{1}{k^2}+\frac{2}{k^3}p^m\pmod{p^{m+1}},
\end{align*}
we have
\begin{align}
\sum_{k=1}^{p^m-1}{'}\frac{1}{k^2}&=\sum_{k=1}^{(p^m-1)/2}{'}\frac{1}{k^2}
+\sum_{k=1}^{(p^m-1)/2}{'}\frac{1}{(p^m-k)^2}\notag\\
&\equiv 2\sum_{k=1}^{(p^m-1)/2}{'}\frac{1}{k^2}+2p^m\sum_{k=1}^{(p^m-1)/2}{'}\frac{1}{k^3}\pmod{p^{m+1}}.
\label{b-11}
\end{align}

By \eqref{b-new-1} and the fact $(p^m-1)/2\equiv (p-1)/2\pmod{p}$, we have
\begin{align}
\sum_{k=1}^{(p^m-1)/2}{'}\frac{1}{k^3}
&\equiv \sum_{k=1}^{(p-1)/2}\frac{1}{k^3}\notag\\
&\equiv -2B_{p-3}\pmod{p},\label{b-12}
\end{align}
where we have used the result \cite[Corollary 5.2]{sunzh-dam-2000} in the last step.

Furthermore, combining \eqref{b-9}, \eqref{b-11} and \eqref{b-12} gives
\begin{align}
\sum_{k=1}^{(p^m-1)/2}{'}\frac{1}{k^2}\equiv \frac{7}{3}p^mB_{p-3}\pmod{p^{m+1}}.\label{b-13}
\end{align}

Thus, by \eqref{b-12} and \eqref{b-13} we have
\begin{align*}
\sum_{\substack{\lfloor k/p^{m}\rfloor=n\\[3pt] \{k/p^{m}\}<p^{m}/2}}{'}\frac{1}{k^2}
&=\sum_{k=1}^{(p^m-1)/2}{'}\frac{1}{(np^m+k)^2}\\
&\equiv \sum_{k=1}^{(p^m-1)/2}{'}\frac{1}{k^2}-2np^m\sum_{k=1}^{(p^m-1)/2}{'}\frac{1}{k^3}\\
&\equiv \frac{12n+7}{3}p^mB_{p-3} \pmod{p^{m+1}},
\end{align*}
and
\begin{align*}
\sum_{\substack{\lfloor k/p^{m}\rfloor=n\\[3pt] \{k/p^{m}\}>p^{m}/2}}{'}\frac{1}{k^2}
&=\sum_{k=1}^{(p^m-1)/2}{'}\frac{1}{((n+1)p^m-k)^2}\\
&\equiv \sum_{k=1}^{(p^m-1)/2}{'}\frac{1}{k^2}+2(n+1)p^m\sum_{k=1}^{(p^m-1)/2}{'}\frac{1}{k^3}\\
&\equiv -\frac{12n+5}{3}p^mB_{p-3} \pmod{p^{m+1}},
\end{align*}
as desired.
\qed

\begin{lem}
Let $p\ge 5$ be a prime and $n,l\in \mathbb{N}$. Then
\begin{align}
&\sum_{\substack{\lfloor k/p^{l+1}\rfloor=n\\[3pt] \{k/p^{l+1}\}<p^{l+1}/2}}{'}
 \frac{1}{k^2}\sum_{j=1}^{\lfloor k/p^l\rfloor}{'}\frac{1}{j}\equiv \frac{1}{3}p^l B_{p-3}\pmod{p^{l+1}},\label{b-14}\\
&\sum_{\substack{\lfloor k/p^{l+1}\rfloor=n\\[3pt] \{k/p^{l+1}\}>p^{l+1}/2}}{'}
 \frac{1}{k^2}\sum_{j=1}^{\lfloor k/p^l\rfloor}{'}\frac{1}{j}\equiv \frac{1}{3}p^l B_{p-3}\pmod{p^{l+1}},\label{b-15}\\
&\sum_{\substack{\lfloor k/p^{l+1}\rfloor=n\\[3pt] \{k/p^{l+1}\}<p^{l+1}/2}}{'}
 \frac{1}{k^2}\sum_{j=1}^{\lfloor 2k/p^l\rfloor}{'}\frac{1}{j}\equiv \frac{4}{3}p^l B_{p-3}\pmod{p^{l+1}},\label{b-16}\\
&\sum_{\substack{\lfloor k/p^{l+1}\rfloor=n\\[3pt] \{k/p^{l+1}\}>p^{l+1}/2}}{'}
 \frac{1}{k^2}\sum_{j=1}^{\lfloor 2k/p^l\rfloor}{'}\frac{1}{j}\equiv \frac{4}{3}p^l B_{p-3} \pmod{p^{l+1}}.\label{b-17}
\end{align}
\end{lem}
{\noindent\it Proof.}
Note that
\begin{align}
\sum_{\substack{\lfloor k/p^{l+1}\rfloor=n\\[3pt] \{k/p^{l+1}\}<p^{l+1}/2}}{'}
 \frac{1}{k^2}\sum_{j=1}^{\lfloor k/p^l\rfloor}{'}\frac{1}{j}
 &=\sum_{k=np^{l+1}+1}^{np^{l+1}+(p^{l+1}-1)/2}{'}
 \frac{1}{k^2}\sum_{j=1}^{\lfloor k/p^l\rfloor}{'}\frac{1}{j}\notag\\[5pt]
 &=\sum_{m=np}^{np+\frac{p-3}{2}}\sum_{\lfloor k/p^l\rfloor=m}{'} \frac{1}{k^2}\sum_{j=1}^{\lfloor k/p^l\rfloor}{'}\frac{1}{j}+\sum_{\substack{{\lfloor k/p^l\rfloor=np+(p-1)/2}\\[3pt] \{k/p^{l}\}<p^{l}/2}}{'} \frac{1}{k^2}\sum_{j=1}^{\lfloor k/p^l\rfloor}{'}\frac{1}{j}.\label{b-18}
\end{align}

By \cite[Theorem 1.2]{sla-camb-2002}, we have
\begin{align*}
\sum_{\lfloor k/p^l\rfloor=m}{'}\frac{1}{k^2}\equiv \frac{2}{3}p^lB_{p-3} \pmod{p^{l+1}}.
\end{align*}
It follows that
\begin{align}
&\sum_{m=np}^{np+\frac{p-3}{2}}\sum_{\lfloor k/p^l\rfloor=m}{'} \frac{1}{k^2}\sum_{j=1}^{\lfloor k/p^l\rfloor}{'}\frac{1}{j}\notag\\[5pt]
&\equiv \frac{2}{3}p^lB_{p-3}\sum_{m=np}^{np+\frac{p-3}{2}}\sum_{j=1}^{m}{'}\frac{1}{j}\notag\\[5pt]
&\equiv \frac{2}{3}p^lB_{p-3}\sum_{k=1}^{(p-3)/2}\sum_{j=1}^k\frac{1}{j} \pmod{p^{l+1}},\label{b-19}
\end{align}
where we have used \eqref{b-new-1} in the last step.

On the other hand, by \eqref{b-new-1} and \eqref{b-7} we have
\begin{align}
&\sum_{\substack{{\lfloor k/p^l\rfloor=np+(p-1)/2}\\[3pt] \{k/p^{l}\}<p^{l}/2}}{'} \frac{1}{k^2}\sum_{j=1}^{\lfloor k/p^l\rfloor}{'}\frac{1}{j}\notag\\[5pt]
&=\sum_{j=1}^{np+(p-1)/2}{'}\frac{1}{j}\sum_{\substack{{\lfloor k/p^l\rfloor=np+(p-1)/2}\\[3pt] \{k/p^{l}\}<p^{l}/2}}{'} \frac{1}{k^2}\notag\\[5pt]
&\equiv \frac{1}{3}p^lB_{p-3}\sum_{j=1}^{np+(p-1)/2}{'}\frac{1}{j}\notag\\[5pt]
&\equiv \frac{1}{3}p^lB_{p-3}\sum_{k=1}^{(p-1)/2}\frac{1}{k}\pmod{p^{l+1}}.\label{b-20}
\end{align}

It follows from \eqref{b-18}, \eqref{b-19} and \eqref{b-20} that
\begin{align}
&\sum_{\substack{\lfloor k/p^{l+1}\rfloor=n\\[3pt] \{k/p^{l+1}\}<p^{l+1}/2}}{'}
 \frac{1}{k^2}\sum_{j=1}^{\lfloor k/p^l\rfloor}{'}\frac{1}{j}\notag\\[5pt]
 &\equiv \frac{2}{3}p^lB_{p-3}\sum_{k=1}^{(p-3)/2}\sum_{j=1}^k\frac{1}{j}+\frac{1}{3}p^lB_{p-3}\sum_{k=1}^{(p-1)/2}\frac{1}{k}\pmod{p^{l+1}}.
 \label{b-21}
\end{align}

Note that
\begin{align}
&2\sum_{k=1}^{(p-3)/2}\sum_{j=1}^k\frac{1}{j}+\sum_{k=1}^{(p-1)/2}\frac{1}{k}\notag\\
&=2\sum_{k=1}^{(p-3)/2}\frac{1}{k}\left(\frac{p-1}{2}-k\right)+\sum_{k=1}^{(p-1)/2}\frac{1}{k}\notag\\
&\equiv -\sum_{k=1}^{(p-3)/2}\frac{1}{k}+\sum_{k=1}^{(p-1)/2}\frac{1}{k}+3\notag\\
&=\frac{2}{p-1}+3\notag\\
&\equiv 1\pmod{p}.\label{b-22}
\end{align}

Then the proof of \eqref{b-14} follows from \eqref{b-21} and \eqref{b-22}.
The proof of \eqref{b-15}--\eqref{b-17} runs analogously, and we omit the details.
\qed

\begin{lem}
Let $p\ge 5$ be a prime and $m\in \mathbb{Z}^{+}$. Then
\begin{align}
\sum_{k=1}^{(p^m-1)/2}{'}\frac{(-1)^{k}}{k^3}\equiv -\frac{1}{4}B_{p-3}\pmod{p}.\label{b-23}
\end{align}
\end{lem}
{\noindent\it Proof.}
By \eqref{b-12}, we have
\begin{align}
\sum_{k=1}^{p^m-1}{'}\frac{(-1)^{k}}{k^3}
&=\frac{1}{4}\sum_{k=1}^{(p^m-1)/2}{'}\frac{1}{k^3}-\sum_{k=1}^{p^m-1}{'}\frac{1}{k^3}\notag\\
&\equiv -\frac{1}{2}B_{p-3}\pmod{p}.\label{b-24}
\end{align}
Note that
\begin{align}
\sum_{k=1}^{p^m-1}{'}\frac{(-1)^{k}}{k^3}&=\sum_{k=1}^{(p^m-1)/2}{'}\frac{(-1)^{k}}{k^3}
+\sum_{k=1}^{(p^m-1)/2}{'}\frac{(-1)^{p^m-k}}{(p^m-k)^3}\notag\\
&\equiv 2\sum_{k=1}^{(p^m-1)/2}{'}\frac{(-1)^{k}}{k^3}\pmod{p}.\label{b-25}
\end{align}
Combining \eqref{b-24} and \eqref{b-25}, we complete the proof of \eqref{b-23}.
\qed

\section{Proof of Theorem \ref{t-1}}
Let
\begin{align*}
S^{(r,s,t)}_n (k)&={n\choose k}^r{n+k\choose k}^s{2k\choose n}^t.
\end{align*}
Note that
\begin{align}
A_{np^m}^{(r,s,t)}&=\sum_{k=0}^{np^m} S^{(r,s,t)}_{np^m}(k)\notag\\[5pt]
&=\sum_{l\ge 1}\sum_{k}{'} S^{(r,s,t)}_{np^m}(kp^l)+\sum_{k=0}^{np^m}{'}S^{(r,s,t)}_{np^m}(k).\label{c-1}
\end{align}

Firstly, we establish a preliminary result on $\sum_{l\ge 1}\sum_{k}{'} S^{(r,s,t)}_{np^m}(kp^l)$ modulo $p^{3m+1}$.

Letting $n\to np^{m-1}$ and $k\to kp^{l-1}$ in \eqref{b-1}, we obtain
\begin{align}
&{np^m\choose kp^l}{\Big/}{np^{m-1}\choose kp^{l-1}}\notag\\[5pt]
&\equiv 1-\frac{nk(np^{m-1}-kp^{l-1})}{3}p^{m+l+1}B_{p-3} \pmod{p^{m+l+\min(m,l)+1}}.
\label{c-2}
\end{align}
Similarly, we have
\begin{align}
&{2kp^l\choose np^m}{\Big /}{2kp^{l-1}\choose np^{m-1}}\notag\\[5pt]
&\equiv 1-\frac{2nk(2kp^{l-1}-np^{m-1})}{3}p^{m+l+1}
B_{p-3}\pmod{p^{m+l+\min(m,l)+1}},\label{c-3}
\end{align}
and
\begin{align}
&{np^m+kp^l\choose kp^l}{\Big /}{np^{m-1}+kp^{l-1}\choose kp^{l-1}}\notag\\[5pt]
&\equiv1-\frac{nk(np^{m-1}+kp^{l-1})}{3}p^{m+l+1}B_{p-3}\pmod{p^{m+l+\min(m,l)+1}}.\label{c-4}
\end{align}
It follows from \eqref{c-2}--\eqref{c-4} that
\begin{align}
 &\frac{S^{(r,s,t)}_{np^m}(kp^l)}{S^{(r,s,t)}_{np^{m-1}}(kp^{l-1})}\notag\\[5pt]
 &\equiv 1-\frac{1}{3}B_{p-3}\left(kn^2p^{2m+l}(r+s-2t)+nk^2p^{m+2l}(s-r+4t)\right)
\pmod{p^{m+l+\min(m,l)+1}}.\label{c-5}
\end{align}

Next, we shall prove that for $p\nmid k$,
\begin{align}
 &S^{(r,s,t)}_{np^m}(kp^l)
 \equiv S^{(r,s,t)}_{np^{m-1}}(kp^{l-1})\notag\\[5pt]
 &-\frac{1}{3}B_{p-3}\left(kn^2p^{2m+l}(r+s-2t)+nk^2p^{m+2l}(s-r+4t)\right)
 S^{(r,s,t)}_{np^{m-1}}(kp^{l-1})\pmod{p^{3m+1}}.\label{c-6}
\end{align}
For $l\ge m$, we have $m+l+\min(m,l)+1\ge 3m+1$, and so \eqref{c-6} follows from \eqref{c-5} directly.
For $l<m$, we have
\begin{align}
m+l+\min(m,l)+1=m+2l+1.\label{c-7}
\end{align}
Noting that for $p\nmid k$ and $r\ge 2$,
\begin{align*}
{np^{m-1}\choose kp^{l-1}}^r=\left(\frac{n}{k}\right)^r{np^{m-1}-1\choose kp^{l-1}-1}^rp^{r(m-l)}
\equiv 0\pmod{p^{2m-2l}},
\end{align*}
we have
\begin{align}
S^{(r,s,t)}_{np^{m-1}}(kp^{l-1})\equiv 0\pmod{p^{2m-2l}}.\label{c-8}
\end{align}
Then \eqref{c-6} follows from \eqref{c-5}, \eqref{c-7} and \eqref{c-8}.

By \eqref{c-6}, we have
\begin{align*}
&\sum_{l\ge 1}\sum_{k}{'} S^{(r,s,t)}_{np^m}(kp^l)\\
&\equiv \sum_{l\ge 1}\sum_{k}{'} S^{(r,s,t)}_{np^{m-1}}(kp^{l-1})\\
&-\frac{1}{3}n^2(r+s-2t)B_{p-3} \sum_{l\ge 1}\sum_{k}{'}p^{2m+l}kS^{(r,s,t)}_{np^{m-1}}(kp^{l-1})\\
&-\frac{1}{3}n(s-r+4t)B_{p-3}\sum_{l\ge 1}\sum_{k}{'}p^{m+2l}k^2S^{(r,s,t)}_{np^{m-1}}(kp^{l-1})
\pmod{p^{3m+1}}.
\end{align*}
Noting that
\begin{align*}
\sum_{l\ge 1}\sum_{k}{'} S^{(r,s,t)}_{np^{m-1}}(kp^{l-1})=A_{np^{m-1}}^{(r,s,t)},
\end{align*}
and $2m+l,m+2l\ge 3m+1$ for $l>m$, we obtain
\begin{align}
&\sum_{l\ge 1}\sum_{k}{'} S^{(r,s,t)}_{np^m}(kp^l)\notag\\
&\equiv A_{np^{m-1}}^{(r,s,t)}\notag\\
&-\frac{1}{3}n^2(r+s-2t)B_{p-3} \sum_{l=1}^{m}\sum_{k}{'}p^{2m+l}kS^{(r,s,t)}_{np^{m-1}}(kp^{l-1})\notag\\
&-\frac{1}{3}n(s-r+4t)B_{p-3}\sum_{l=1}^{m}\sum_{k}{'}p^{m+2l}k^2S^{(r,s,t)}_{np^{m-1}}(kp^{l-1})
\pmod{p^{3m+1}}.\label{c-9}
\end{align}

For $p\nmid k$, $r\ge 2$ and $l<m$, we have
\begin{align*}
p^{2m+l}{np^{m-1}\choose kp^{l-1}}^r=p^{(r+2)m+(1-r)l}\left(\frac{n}{k}\right)^r{np^{m-1}-1\choose kp^{l-1}-1}^r
\equiv 0\pmod{p^{3m+1}},
\end{align*}
and so
\begin{align*}
\sum_{k}{'}p^{2m+l}kS^{(r,s,t)}_{np^{m-1}}(kp^{l-1})
\equiv 0\pmod{p^{3m+1}}.
\end{align*}
It follows that
\begin{align}
&\sum_{l=1}^{m}\sum_{k}{'}p^{2m+l} kS^{(r,s,t)}_{np^{m-1}}(kp^{l-1})\notag\\
&\equiv p^{3m} \sum_{k}{'}kS^{(r,s,t)}_{np^{m-1}}(kp^{m-1})\notag\\
&\equiv p^{3m} \sum_{k=0}^n {'}k S^{(r,s,t)}_n (k)\notag\\
&\equiv p^{3m} \sum_{k=0}^n k S^{(r,s,t)}_n (k) \pmod{p^{3m+1}},\label{c-10}
\end{align}
where we have used the modulo $p$ version of \eqref{b-1}.

Combining \eqref{c-9} and \eqref{c-10} gives
\begin{align}
&\sum_{l\ge 1}\sum_{k}{'} S^{(r,s,t)}_{np^m}(kp^l)\notag\\
&\equiv A_{np^{m-1}}^{(r,s,t)}\notag\\
&-\frac{1}{3}p^{3m}n^2(r+s-2t)B_{p-3}  \sum_{k=0}^n k S^{(r,s,t)}_n (k)\notag\\
&-\frac{1}{3}n(s-r+4t)B_{p-3}\sum_{l=1}^{m}\sum_{k}{'}p^{m+2l}k^2S^{(r,s,t)}_{np^{m-1}}(kp^{l-1})
\pmod{p^{3m+1}}.\label{c-11}
\end{align}

Secondly, we rewrite $\sum_{k}{'}S^{(r,s,t)}_{np^m}(k)$ as
\begin{align}
\sum_{k=0}^{np^m}{'}S^{(r,s,t)}_{np^m}(k)
=n^rp^{rm}\sum_{k}{'} \frac{1}{k^r}{np^m-1\choose k-1}^r{np^m+k\choose k}^s{2k\choose np^m}^t. \label{c-12}
\end{align}

Let $Y^{(r,s,t)}_{m,n,p}$ denote the last double sum on the right-hand side of \eqref{c-11} and
$Z^{(r,s,t)}_{m,n,p}$ denote the right-hand side of \eqref{c-12}, namely,
\begin{align*}
&Y^{(r,s,t)}_{m,n,p}=\sum_{l=1}^{m}\sum_{k}{'}p^{m+2l}k^2S^{(r,s,t)}_{np^{m-1}}(kp^{l-1}),\\[5pt]
&Z^{(r,s,t)}_{m,n,p}=n^rp^{rm}\sum_{k}{'} \frac{1}{k^r}{np^m-1\choose k-1}^r{np^m+k\choose k}^s{2k\choose np^m}^t.
\end{align*}
It follows from \eqref{c-1}, \eqref{c-11} and \eqref{c-12} that
\begin{align}
A_{np^m}^{(r,s,t)}
&\equiv A_{np^{m-1}}^{(r,s,t)}\notag\\
&-\frac{1}{3}p^{3m}n^2(r+s-2t)B_{p-3}  \sum_{k=0}^n k S^{(r,s,t)}_n (k)\notag\\
&-\frac{1}{3}n(s-r+4t)B_{p-3}Y^{(r,s,t)}_{m,n,p}+Z^{(r,s,t)}_{m,n,p}\pmod{p^{3m+1}}.\label{c-new-1}
\end{align}
In order to establish the Theorem \ref{t-1}, it suffices to determine $Y^{(r,s,t)}_{m,n,p}$ and $Z^{(r,s,t)}_{m,n,p}$ modulo $p^{3m+1}$. Substitution of the computed results $Y^{(r,s,t)}_{m,n,p}$ and $Z^{(r,s,t)}_{m,n,p}$ modulo $p^{3m+1}$ into \eqref{c-new-1} sets, after simplifications, up the Theorem \ref{t-1}. The rest of the proof is divided into two parts.

\subsection{$Y^{(r,s,t)}_{m,n,p}$ modulo $p^{3m+1}$}
We shall distinguish three cases to determine $Y^{(r,s,t)}_{m,n,p}$ modulo $p^{3m+1}$.

{\noindent\bf Case 1 $m=1$.}

It is easy to see that
\begin{align*}
Y^{(r,s,t)}_{1,n,p}&=p^{3}\sum_{k=0}^n{'}k^2S^{(r,s,t)}_{n}(k)\\
&\equiv p^{3}\sum_{k=0}^nk^2S^{(r,s,t)}_{n}(k)\pmod{p^4}.
\end{align*}

{\noindent\bf Case 2 $m\ge 2$ and $r=2$.}

Note that
\begin{align}
Y^{(2,s,t)}_{m,n,p}
&= p^{3m}n^2\sum_{l=1}^{m}\sum_{k}{'} {np^{m-1}-1\choose kp^{l-1}-1}^2{np^{m-1}+kp^{l-1}\choose kp^{l-1}}^s{2kp^{l-1}\choose np^{m-1}}^t\notag\\[5pt]
&=p^{3m}n^2\sum_{k}{'} {np^{m-1}-1\choose kp^{m-1}-1}^2{np^{m-1}+kp^{m-1}\choose kp^{m-1}}^s{2kp^{m-1}\choose np^{m-1}}^t\notag\\[5pt]
&+p^{3m}n^2\sum_{l=1}^{m-1}\sum_{k}{'} {np^{m-1}-1\choose kp^{l-1}-1}^2{np^{m-1}+kp^{l-1}\choose kp^{l-1}}^s{2kp^{l-1}\choose np^{m-1}}^t.\label{c-13}
\end{align}

By a repeated use of the modulo $p$ version of \eqref{b-2}, we obtain
\begin{align}
&p^{3m}n^2\sum_{k}{'} {np^{m-1}-1\choose kp^{m-1}-1}^2{np^{m-1}+kp^{m-1}\choose kp^{m-1}}^s{2kp^{m-1}\choose np^{m-1}}^t\notag\\[5pt]
&\equiv p^{3m}n^2\sum_{k}{'} {n-1\choose k-1}^2{n+k\choose k}^s{2k\choose n}^t\notag\\[5pt]
&=p^{3m}\sum_{k}{'} k^2S^{(2,s,t)}_{n}(k)\notag\\[5pt]
&\equiv p^{3m}\sum_{k=0}^n k^2S^{(2,s,t)}_{n}(k) \pmod{p^{3m+1}}.
\label{c-14}
\end{align}

On the other hand, we have
\begin{align}
&p^{3m}n^2\sum_{l=1}^{m-1}\sum_{k}{'} {np^{m-1}-1\choose kp^{l-1}-1}^2{np^{m-1}+kp^{l-1}\choose kp^{l-1}}^s{2kp^{l-1}\choose np^{m-1}}^t\notag\\[5pt]
&=p^{3m}n^2\sum_{p^{m-1}\nmid k}{np^{m-1}-1\choose k-1}^2{np^{m-1}+k\choose k}^s{2k\choose np^{m-1}}^t\notag\\[5pt]
&\equiv p^{3m}n^2\sum_{p^{m-1}\nmid k}{n-1\choose \lfloor k/p^{m-1}\rfloor }^2{n+\lfloor k/p^{m-1}\rfloor\choose \lfloor k/p^{m-1}\rfloor}^s{\lfloor 2k/p^{m-1}\rfloor\choose n}^t\pmod{p^{3m+1}},
\label{c-15}
\end{align}
where we have used the modulo $p$ version of \eqref{b-2}. Note that for $m\in \mathbb{Z}^{+}$,
\begin{align}
\lfloor 2k/p^{m}\rfloor=
\begin{cases}
2\lfloor k/p^{m}\rfloor\quad &\text{if $\{k/p^{m}\}<p^{m}/2$},\\[5pt]
2\lfloor k/p^{m}\rfloor+1 \quad &\text{if $\{k/p^{m}\}>p^{m}/2$}.\label{c-16}
\end{cases}
\end{align}
It follows from \eqref{c-15} and \eqref{c-16} that
\begin{align}
&p^{3m}n^2\sum_{l=1}^{m-1}\sum_{k}{'} {np^{m-1}-1\choose kp^{l-1}-1}^2{np^{m-1}+kp^{l-1}\choose kp^{l-1}}^s{2kp^{l-1}\choose np^{m-1}}^t\notag\\[5pt]
&\equiv p^{3m}n^2\sum_{N}\sum_{\substack{p^{m-1}\nmid k\\[3pt] \lfloor k/p^{m-1}\rfloor=N\\[3pt] \{k/p^{m-1}\}<p^{m-1}/2}}{n-1\choose N }^2{n+N\choose N}^s{2N\choose n}^t\notag\\[5pt]
&+p^{3m}n^2\sum_{N}\sum_{\substack{p^{m-1}\nmid k\\[3pt] \lfloor k/p^{m-1}\rfloor=N\\[3pt] \{k/p^{m-1}\}>p^{m-1}/2}}{n-1\choose N }^2{n+N\choose N}^s{2N+1\choose n}^t\notag\\[5pt]
&\equiv-\frac{1}{2} p^{3m}n^2\sum_{k=0}^{n-1}{n-1\choose k }^2{n+k\choose k}^s{2k\choose n}^t\notag\\[5pt]
&-\frac{1}{2} p^{3m}n^2\sum_{k=0}^{n-1}{n-1\choose k }^2{n+k\choose k}^s{2k+1\choose n}^t\pmod{p^{3m+1}},
\label{c-17}
\end{align}
where we have used the following fact:
\begin{align*}
\sum_{\substack{p^{m-1}\nmid k\\[3pt] \lfloor k/p^{m-1}\rfloor=N\\[3pt] \{k/p^{m-1}\}<p^{m-1}/2}}1
=\sum_{\substack{p^{m-1}\nmid k\\[3pt] \lfloor k/p^{m-1}\rfloor=N\\[3pt] \{k/p^{m-1}\}>p^{m-1}/2}}1
=\frac{p^{m-1}-1}{2}\equiv -\frac{1}{2}\pmod{p}.
\end{align*}

It follows from \eqref{c-13}, \eqref{c-14} and \eqref{c-17} that
\begin{align*}
Y^{(2,s,t)}_{m,n,p}
&\equiv p^{3m}\sum_{k=0}^n k^2S^{(2,s,t)}_{n}(k)\\
&-\frac{1}{2} p^{3m}n^2\sum_{k=0}^{n-1}{n-1\choose k }^2{n+k\choose k}^s{2k\choose n}^t\\
&-\frac{1}{2} p^{3m}n^2\sum_{k=0}^{n-1}{n-1\choose k }^2{n+k\choose k}^s{2k+1\choose n}^t\pmod{p^{3m+1}}.
\end{align*}

{\noindent\bf Case 3 $m\ge 2$ and $r\ge 3$.}

For $p\nmid k$ and $l<m$, we have
\begin{align*}
p^{m+2l}k^2 {np^{m-1}\choose kp^{l-1}}^r
=p^{(r+1)m+(2-r)l}k^2 \left(\frac{n}{k}\right)^r{np^{m-1}-1\choose kp^{l-1}-1}^r
\equiv 0\pmod{p^{3m+1}},
\end{align*}
and so
\begin{align*}
p^{m+2l}k^2S^{(r,s,t)}_{np^{m-1}}(kp^{l-1})\equiv 0\pmod{p^{3m+1}}.
\end{align*}
It follows that
\begin{align*}
Y^{(r,s,t)}_{m,n,p}
&\equiv p^{3m}\sum_{k}{'}k^2S^{(r,s,t)}_{np^{m-1}}(kp^{m-1})\\
&\equiv p^{3m}\sum_{k=0}^n{'}k^2S^{(r,s,t)}_{n}(k)\\
&\equiv p^{3m}\sum_{k=0}^n k^2S^{(r,s,t)}_{n}(k)\pmod{p^{3m+1}},
\end{align*}
where we have used the modulo $p$ version of \eqref{b-1}.

\subsection{$Z^{(r,s,t)}_{m,n,p}$ modulo $p^{3m+1}$}
We shall distinguish four cases to determine $Z^{(r,s,t)}_{m,n,p}$ modulo $p^{3m+1}$.

{\noindent\bf Case 1 $r\ge 4$.}

In this case, we have $p^{rm}\equiv 0\pmod{p^{3m+1}}$, and so
\begin{align*}
Z^{(r,s,t)}_{m,n,p}\equiv 0\pmod{p^{3m+1}}.
\end{align*}

{\noindent\bf Case 2 $r=3$.}

By a repeated use of the modulo $p$ version of \eqref{b-2}, we obtain
\begin{align}
&Z^{(3,s,t)}_{m,n,p}\notag\\[5pt]
&\equiv -n^3p^{3m}\sum_{k}{'} \frac{(-1)^{k+\lfloor k/p^m\rfloor}}{k^3}{n-1\choose \lfloor k/p^m\rfloor}^3{n+\lfloor k/p^m\rfloor\choose \lfloor k/p^m\rfloor}^s{\lfloor 2k/p^m\rfloor\choose n}^t\notag\\[5pt]
&=-n^3p^{3m}\sum_{N}{n-1\choose N}^3{n+N\choose N}^s{2N\choose n}^t \sum_{\substack{\lfloor k/p^m\rfloor=N\\[3pt]\{k/p^m\}<p^m/2 }}{'}\frac{(-1)^{N+k}}{k^3}\notag\\[5pt]
&-n^3p^{3m}\sum_{N}{n-1\choose N}^3{n+N\choose N}^s{2N+1\choose n}^t \sum_{\substack{\lfloor k/p^m\rfloor=N\\[3pt]\{k/p^m\}>p^m/2 }}{'}\frac{(-1)^{N+k}}{k^3}\pmod{p^{3m+1}},\label{c-new-2}
\end{align}
where we have used the fact that $\lfloor (k-1)/p^m\rfloor=\lfloor k/p^m\rfloor$ for $p\nmid k$.

By \eqref{b-23}, we have
\begin{align}
\sum_{\substack{\lfloor k/p^m\rfloor=N\\[3pt]\{k/p^m\}<p^m/2 }}{'}\frac{(-1)^{N+k}}{k^3}
=&\sum_{k=1}^{(p^m-1)/2}{'}\frac{(-1)^{N+Np^m+k}}{(Np^m+k)^3}\notag\\[5pt]
&\equiv \sum_{k=1}^{(p^m-1)/2}{'}\frac{(-1)^{k}}{k^3}\notag\\[5pt]
&\equiv -\frac{1}{4}B_{p-3}\pmod{p},\label{c-new-3}
\end{align}
and
\begin{align}
\sum_{\substack{\lfloor k/p^m\rfloor=N\\[3pt]\{k/p^m\}>p^m/2 }}{'}\frac{(-1)^{N+k}}{k^3}
=&\sum_{k=1}^{(p^m-1)/2}{'}\frac{(-1)^{N+(N+1)p^m-k}}{((N+1)p^m-k)^3}\notag\\[5pt]
&\equiv \sum_{k=1}^{(p^m-1)/2}{'}\frac{(-1)^{k}}{k^3}\notag\\[5pt]
&\equiv -\frac{1}{4}B_{p-3}\pmod{p}.\label{c-new-4}
\end{align}
It follows from \eqref{c-new-2}--\eqref{c-new-4} that
\begin{align*}
&Z^{(3,s,t)}_{m,n,p}\notag\\[5pt]
&\equiv \frac{1}{4}p^{3m}B_{p-3}\sum_{k}{n-1\choose k}^3{n+k\choose k}^s{2k\choose n}^t n^3 \notag\\[5pt]
&+\frac{1}{4}p^{3m}B_{p-3}\sum_{k}{n-1\choose k}^3{n+k\choose k}^s{2k+1\choose n}^tn^3\pmod{p^{3m+1}}.
\end{align*}

{\noindent\bf Case 3 $r=2$ and $m=1$.}

By \eqref{b-2}, we have
\begin{align}
Z^{(2,s,t)}_{1,n,p}&=n^2p^{2}\sum_{k}{'} \frac{1}{k^2}{np-1\choose k-1}^2{np+k\choose k}^s{2k\choose n}^t\notag\\[5pt]
&\equiv n^2p^{2}\sum_{k}{'} \frac{1}{k^2}{n-1\choose \lfloor k/p\rfloor}^2{n+\lfloor k/p\rfloor\choose \lfloor k/p\rfloor}^s{\lfloor 2k/p\rfloor\choose n}^t\notag\\[5pt]
&\times \left(1-2np\sum_{j=1}^{k-1}{'}\frac{1}{j}
+snp\sum_{j=1}^k{'}\frac{1}{j}+tnp\sum_{j=1}^{2k}{'}\frac{1}{j}\right)\pmod{p^{4}},\label{c-19}
\end{align}
where we have used the fact that $\lfloor (k-1)/p\rfloor=\lfloor k/p\rfloor$ for $p\nmid k$.

By \eqref{b-7}, \eqref{b-8} and \eqref{c-16}, we have
\begin{align}
&\sum_{k}{'} \frac{1}{k^2}{n-1\choose \lfloor k/p\rfloor}^2{n+\lfloor k/p\rfloor\choose \lfloor k/p\rfloor}^s{\lfloor 2k/p\rfloor\choose n}^t\notag\\[5pt]
&=\sum_{N}{n-1\choose N}^2{n+N\choose N}^s{2N\choose n}^t\sum_{\substack{\lfloor k/p\rfloor=N\\[3pt]  \{k/p\}<p/2}}{'}\frac{1}{k^2}\notag\\[5pt]
&+\sum_{N}{n-1\choose N}^2{n+N\choose N}^s{2N+1\choose n}^t\sum_{\substack{\lfloor k/p\rfloor=N\\[3pt]  \{k/p\}>p/2}}{'}\frac{1}{k^2}\notag\\[5pt]
&\equiv \frac{1}{3}pB_{p-3}\sum_{k}{n-1\choose k}^2{n+k\choose k}^s{2k\choose n}^t(12k+7)\notag\\[5pt]
&-\frac{1}{3}pB_{p-3}\sum_{k}{n-1\choose k}^2{n+k\choose k}^s{2k+1\choose n}^t(12k+5)\pmod{p^2}.
\label{c-20}
\end{align}

By \eqref{c-16}, we have
\begin{align}
&\sum_{k}{'} \frac{1}{k^2}{n-1\choose \lfloor k/p\rfloor}^2{n+\lfloor k/p\rfloor\choose \lfloor k/p\rfloor}^s{\lfloor 2k/p\rfloor\choose n}^t
\left(-2\sum_{j=1}^{k-1}{'}\frac{1}{j}
+s\sum_{j=1}^k{'}\frac{1}{j}+t\sum_{j=1}^{2k}{'}\frac{1}{j}\right)\notag\\[5pt]
&=\sum_{k}{'} \frac{1}{k^2}{n-1\choose \lfloor k/p\rfloor}^2{n+\lfloor k/p\rfloor\choose \lfloor k/p\rfloor}^s{\lfloor 2k/p\rfloor\choose n}^t
\left(\frac{2}{k}+(s-2)\sum_{j=1}^k{'}\frac{1}{j}+t\sum_{j=1}^{2k}{'}\frac{1}{j}\right)\notag\\[5pt]
&=\sum_{N} {n-1\choose N}^2{n+N\choose N}^s{2N\choose n}^t\sum_{\substack{\lfloor k/p\rfloor=N\\[3pt]  \{k/p\}<p/2}}{'}\frac{1}{k^2} \left(\frac{2}{k}+(s-2)\sum_{j=1}^k{'}\frac{1}{j}+t\sum_{j=1}^{2k}{'}\frac{1}{j}\right)\notag\\[5pt]
&+\sum_{N} {n-1\choose N}^2{n+N\choose N}^s{2N+1\choose n}^t\sum_{\substack{\lfloor k/p\rfloor=N\\[3pt]  \{k/p\}>p/2}}{'}\frac{1}{k^2} \left(\frac{2}{k}+(s-2)\sum_{j=1}^k{'}\frac{1}{j}+t\sum_{j=1}^{2k}{'}\frac{1}{j}\right).\label{c-21}
\end{align}
By the result \cite[Corollary 5.2]{sunzh-dam-2000} and the fact $\sum_{k=1}^{p-1}\frac{1}{k^3}\equiv 0\pmod{p}$, we have
\begin{align}
&\sum_{\substack{\lfloor k/p\rfloor=N\\[3pt]  \{k/p\}<p/2}}{'}\frac{1}{k^3}\equiv -2B_{p-3}\pmod{p},\label{c-22}\\[5pt]
&\sum_{\substack{\lfloor k/p\rfloor=N\\[3pt]  \{k/p\}>p/2}}{'}\frac{1}{k^3}\equiv 2B_{p-3}\pmod{p}.\label{c-23}
\end{align}
It follows from \eqref{b-14}--\eqref{b-17} and \eqref{c-21}--\eqref{c-23} that
\begin{align}
&\sum_{k}{'} \frac{1}{k^2}{n-1\choose \lfloor k/p\rfloor}^2{n+\lfloor k/p\rfloor\choose \lfloor k/p\rfloor}^s{\lfloor 2k/p\rfloor\choose n}^t
\left(-2\sum_{j=1}^{k-1}{'}\frac{1}{j}
+s\sum_{j=1}^k{'}\frac{1}{j}+t\sum_{j=1}^{2k}{'}\frac{1}{j}\right)\notag\\[5pt]
&\equiv \frac{1}{3}B_{p-3}\sum_{k} {n-1\choose k}^2{n+k\choose k}^s{2k\choose n}^t(4t+s-14)\notag\\[5pt]
&+\frac{1}{3}B_{p-3}\sum_{k} {n-1\choose k}^2{n+k\choose k}^s{2k+1\choose n}^t(4t+s+10)\pmod{p}.\label{c-24}
\end{align}

Finally, combining \eqref{c-19}, \eqref{c-20} and \eqref{c-24} gives
\begin{align*}
Z^{(2,s,t)}_{1,n,p}\equiv
&\frac{1}{3}p^3B_{p-3}\sum_{k}{n-1\choose k}^2{n+k\choose k}^s{2k\choose n}^t n^2(12k+7)\notag\\[5pt]
&-\frac{1}{3}p^3B_{p-3}\sum_{k}{n-1\choose k}^2{n+k\choose k}^s{2k+1\choose n}^t n^2(12k+5)\\[5pt]
&+\frac{1}{3}p^3B_{p-3}\sum_{k} {n-1\choose k}^2{n+k\choose k}^s{2k\choose n}^tn^3(4t+s-14)\notag\\[5pt]
&+\frac{1}{3}p^3B_{p-3}\sum_{k} {n-1\choose k}^2{n+k\choose k}^s{2k+1\choose n}^t n^3(4t+s+10)\pmod{p^4}.
\end{align*}

{\noindent\bf Case 4 $r=2$ and $m\ge 2$.}

By \eqref{b-2}, we have
\begin{align}
Z^{(2,s,t)}_{m,n,p}
&\equiv n^2p^{2m}\sum_{k}{'} \frac{1}{k^2}{np^{m-1}-1\choose \lfloor k/p\rfloor}^2
{np^{m-1}+\lfloor k/p\rfloor\choose \lfloor k/p\rfloor}^s
{\lfloor 2k/p\rfloor \choose np^{m-1}}^t\notag\\[5pt]
&\times\left(1-2np^m\sum_{j=1}^{k-1}{'}\frac{1}{j}
+snp^m\sum_{j=1}^k{'}\frac{1}{j}+tnp^m\sum_{j=1}^{2k}{'}\frac{1}{j}\right)\pmod{p^{3m+1}},\label{c-25}
\end{align}
where we have used the fact that $\lfloor (k-1)/p\rfloor=\lfloor k/p\rfloor$ for $p\nmid k$.

Let
\begin{align*}
T_{s,t}(k)=s\sum_{j=1}^k{'}\frac{1}{j}+t\sum_{j=1}^{2k}{'}\frac{1}{j}-2\sum_{j=1}^{k-1}{'}\frac{1}{j}.
\end{align*}
Firstly, we shall determine the following sum modulo $p$:
\begin{align*}
{\sum_{k}}' \frac{1}{k^2}{np^{m-1}-1\choose \lfloor k/p\rfloor}^2
{np^{m-1}+\lfloor k/p\rfloor\choose \lfloor k/p\rfloor}^s
{\lfloor 2k/p\rfloor \choose np^{m-1}}^tT_{s,t}(k).
\end{align*}

Note that
\begin{align}
&{\sum_{k}}' \frac{1}{k^2}{np^{m-1}-1\choose \lfloor k/p\rfloor}^2
{np^{m-1}+\lfloor k/p\rfloor\choose \lfloor k/p\rfloor}^s
{\lfloor 2k/p\rfloor \choose np^{m-1}}^t T_{s,t}(k) \notag\\[5pt]
&=\sum_{N}{np^{m-1}-1\choose N}^2
{np^{m-1}+N\choose N}^s
{2N \choose np^{m-1}}^t {\sum_{\substack{\lfloor k/p\rfloor=N\\[3pt] \{k/p\}<p/2}}}' \frac{1}{k^2}T_{s,t}(k)\notag\\[5pt]
&+\sum_{N}{np^{m-1}-1\choose N}^2
{np^{m-1}+N\choose N}^s
{2N+1 \choose np^{m-1}}^t {\sum_{\substack{\lfloor k/p\rfloor=N\\[3pt] \{k/p\}>p/2}}}' \frac{1}{k^2}T_{s,t}(k).\label{c-26}
\end{align}
By using \eqref{b-14}--\eqref{b-17}, \eqref{c-22} and \eqref{c-23}, we obtain
\begin{align}
&\sum_{\substack{\lfloor k/p\rfloor=N\\[3pt] \{k/p\}<p/2}}{'} \frac{1}{k^2}T_{s,t}(k)
\equiv \frac{1}{2}(3t-s-6)B_{p-3}\pmod{p},\label{c-27}\\[5pt]
&\sum_{\substack{\lfloor k/p\rfloor=N\\[3pt] \{k/p\}>p/2}}{'} \frac{1}{k^2}T_{s,t}(k)
\equiv \frac{1}{2}(3s+5t+2)B_{p-3}\pmod{p}.\label{c-28}
\end{align}
By a repeated use of the modulo $p$ version of \eqref{b-2}, we obtain
\begin{align}
&\sum_{N}{np^{m-1}-1\choose N}^2
{np^{m-1}+N\choose N}^s
{2N \choose np^{m-1}}^t\notag\\[5pt]
&\equiv \sum_{N}{n-1\choose \lfloor N/p^{m-1}\rfloor}^2
{n+\lfloor N/p^{m-1}\rfloor\choose \lfloor N/p^{m-1}\rfloor}^s
{\lfloor 2N/p^{m-1}\rfloor \choose n}^t\notag\\[5pt]
&=\sum_{M} \sum_{\substack{\lfloor N/p^{m-1}\rfloor=M\\[3pt] \{N/p^{m-1}\}<p^{m-1}/2}}{n-1\choose M}^2
{n+M\choose M}^s
{2M \choose n}^t\notag\\[5pt]
&+\sum_{M} \sum_{\substack{\lfloor N/p^{m-1}\rfloor=M\\[3pt] \{N/p^{m-1}\}>p^{m-1}/2}}{n-1\choose M}^2
{n+M\choose M}^s
{2M+1 \choose n}^t\notag\\[5pt]
&\equiv \frac{1}{2}\sum_{k}{n-1\choose k}^2
{n+k\choose k}^s
{2k \choose n}^t
-\frac{1}{2}\sum_{k}{n-1\choose k}^2
{n+k\choose k}^s
{2k+1 \choose n}^t\pmod{p},\label{c-29}
\end{align}
where we have used the fact that
\begin{align*}
&\sum_{\substack{\lfloor N/p^{m-1}\rfloor=M\\[3pt] \{N/p^{m-1}\}<p^{m-1}/2}}1=\frac{p^{m-1}+1}{2}\equiv \frac{1}{2}\pmod{p},\\[5pt]
&\sum_{\substack{\lfloor N/p^{m-1}\rfloor=M\\[3pt] \{N/p^{m-1}\}>p^{m-1}/2}}1=\frac{p^{m-1}-1}{2}\equiv -\frac{1}{2}\pmod{p}.
\end{align*}
Similarly, we have
\begin{align}
&\sum_{N}{np^{m-1}-1\choose N}^2
{np^{m-1}+N\choose N}^s
{2N+1 \choose np^{m-1}}^t\notag\\[5pt]
&\equiv \frac{1}{2}\sum_{k}{n-1\choose k}^2
{n+k\choose k}^s
{2k+1 \choose n}^t
-\frac{1}{2}\sum_{k}{n-1\choose k}^2
{n+k\choose k}^s
{2k \choose n}^t\pmod{p}.\label{c-30}
\end{align}
It follows from \eqref{c-26}--\eqref{c-30} that
\begin{align}
&{\sum_{k}}' \frac{1}{k^2}{np^{m-1}-1\choose \lfloor k/p\rfloor}^2
{np^{m-1}+\lfloor k/p\rfloor\choose \lfloor k/p\rfloor}^s
{\lfloor 2k/p\rfloor \choose np^{m-1}}^tT_{s,t}(k)\notag\\[5pt]
&\equiv \frac{1}{2}B_{p-3}\sum_{k}{n-1\choose k}^2
{n+k\choose k}^s{2k+1 \choose n}^t (t+2s+4) \notag\\[5pt]
&-\frac{1}{2}B_{p-3}\sum_{k}{n-1\choose k}^2
{n+k\choose k}^s{2k \choose n}^t(t+2s+4)\pmod{p}.\label{c-31}
\end{align}

For $l<m$, let
\begin{align*}
C^{(2,s,t)}_{n,p,m,l}(k)={np^{m-l}-1\choose \lfloor k/p^l\rfloor}^2
{np^{m-l}+\lfloor k/p^l\rfloor\choose \lfloor k/p^l\rfloor}^s
{\lfloor 2k/p^l\rfloor \choose np^{m-l}}^t.
\end{align*}
Secondly, we shall determine the following sum modulo $p^{m+1}$:
\begin{align*}
\sum_{k}{'} \frac{1}{k^2}{np^{m-1}-1\choose \lfloor k/p\rfloor}^2
{np^{m-1}+\lfloor k/p\rfloor\choose \lfloor k/p\rfloor}^s
{\lfloor 2k/p\rfloor \choose np^{m-1}}^t=\sum_{k}{'} \frac{1}{k^2}C^{(2,s,t)}_{n,p,m,1}(k).
\end{align*}

By \eqref{b-2}, we have
\begin{align}
&C^{(2,s,t)}_{n,p,m,l}(k)\equiv C^{(2,s,t)}_{n,p,m,l+1}(k)\notag\\[5pt]
&\times\left(1-2np^{m-l}\sum_{j=1}^{\lfloor k/p^l\rfloor}{'}\frac{1}{j}
+snp^{m-l}\sum_{j=1}^{\lfloor k/p^l\rfloor}{'}\frac{1}{j}+tnp^{m-l}\sum_{j=1}^{\lfloor 2k/p^l\rfloor}{'}\frac{1}{j}\right)\pmod{p^{m-l+1}}.\label{c-32}
\end{align}
By \eqref{b-7} and \eqref{b-8}, we have
\begin{align}
{\sum_{\substack{\lfloor k/p^l\rfloor=N\\[3pt] \{k/p^l\}<p^l/2}}}'
 \frac{1}{k^2}\equiv {\sum_{\substack{\lfloor k/p^l\rfloor=N\\[3pt] \{k/p^l\}>p^l/2}}}'
 \frac{1}{k^2}\equiv 0 \pmod{p^l}.\label{c-33}
\end{align}
Moreover, we have
\begin{align}
C^{(2,s,t)}_{n,p,m,l}(k)=
\begin{cases}
{np^{m-l}-1\choose N}^2{np^{m-l}+N\choose N}^s{2N \choose np^{m-l}}^t&\quad\text{if $\lfloor k/p^l\rfloor=N$ and $\{k/p^l\}<p^l/2$,}\\[10pt]
{np^{m-l}-1\choose N}^2{np^{m-l}+N\choose N}^s{2N+1 \choose np^{m-l}}^t&\quad\text{if $\lfloor k/p^l\rfloor=N$ and $\{k/p^l\}>p^l/2$.}\label{c-34}
\end{cases}
\end{align}
Note that
\begin{align}
{\sum_{k}}' \frac{1}{k^2}C^{(2,s,t)}_{n,p,m,l}(k)
=\sum_{N}{\sum_{\substack{\lfloor k/p^l\rfloor=N\\[3pt] \{k/p^l\}<p^l/2}}}{'} \frac{1}{k^2}C^{(2,s,t)}_{n,p,m,l}(k)+\sum_{N}{\sum_{\substack{\lfloor k/p^l\rfloor=N\\[3pt] \{k/p^l\}>p^l/2}}}{'} \frac{1}{k^2}C^{(2,s,t)}_{n,p,m,l}(k).\label{c-35}
\end{align}
Applying \eqref{c-32}--\eqref{c-34} to the right-hand side of \eqref{c-35}, we arrive at
\begin{align}
&{\sum_{k}}' \frac{1}{k^2}C^{(2,s,t)}_{n,p,m,l}(k)
\equiv {\sum_{k}}'\frac{1}{k^2}C^{(2,s,t)}_{n,p,m,l+1}(k)\notag\\[5pt]
&\times \left(1-2np^{m-l}\sum_{j=1}^{\lfloor k/p^l\rfloor}{'}\frac{1}{j}
+snp^{m-l}\sum_{j=1}^{\lfloor k/p^l\rfloor}{'}\frac{1}{j}+tnp^{m-l}\sum_{j=1}^{\lfloor 2k/p^l\rfloor}{'}\frac{1}{j}\right)\pmod{p^{m+1}}.\label{c-36}
\end{align}
Note that
\begin{align}
&{\sum_{k}}' \frac{1}{k^2}C^{(2,s,t)}_{n,p,m,l+1}(k)\sum_{j=1}^{\lfloor k/p^l\rfloor}{'}\frac{1}{j}\notag\\[5pt]
&=\sum_{N}{np^{m-l-1}-1\choose N}^2
{np^{m-l-1}+N\choose N}^s{2N \choose np^{m-l-1}}^t{\sum_{\substack{\lfloor k/p^{l+1}\rfloor=N\\[3pt] \{k/p^{l+1}\}<p^{l+1}/2}}}'
 \frac{1}{k^2}\sum_{j=1}^{\lfloor k/p^l\rfloor}{'}\frac{1}{j}\notag\\[5pt]
&+\sum_{N}{np^{m-l-1}-1\choose N}^2
{np^{m-l-1}+N\choose N}^s{2N+1 \choose np^{m-l-1}}^t{\sum_{\substack{\lfloor k/p^{l+1}\rfloor=N\\[3pt] \{k/p^{l+1}\}>p^{l+1}/2}}}'
 \frac{1}{k^2}\sum_{j=1}^{\lfloor k/p^l\rfloor}{'}\frac{1}{j}.\label{c-37}
\end{align}
By \eqref{c-29} and \eqref{c-30}, for $m>l+1$ we have
\begin{align}
&\sum_{N}{np^{m-l-1}-1\choose N}^2{np^{m-l-1}+N\choose N}^s{2N \choose np^{m-l-1}}^t\notag\\[5pt]
&+\sum_{N}{np^{m-l-1}-1\choose N}^2
{np^{m-l-1}+N\choose N}^s{2N+1 \choose np^{m-l-1}}^t\equiv 0\pmod{p}.\label{c-38}
\end{align}
Applying \eqref{b-14}, \eqref{b-15} and \eqref{c-38} to the right-hand side of \eqref{c-37}, we deduce that for $m>l+1$,
\begin{align}
{\sum_{k}}' \frac{1}{k^2}C^{(2,s,t)}_{n,p,m,l+1}(k)\sum_{j=1}^{\lfloor k/p^l\rfloor}{'}\frac{1}{j}\equiv 0\pmod{p^{l+1}}.\label{c-39}
\end{align}
Similarly, from \eqref{b-16}, \eqref{b-17} and \eqref{c-38}, we derive that for $m>l+1$,
\begin{align}
{\sum_{k}}' \frac{1}{k^2}C^{(2,s,t)}_{n,p,m,l+1}(k)\sum_{j=1}^{\lfloor 2k/p^l\rfloor}{'}\frac{1}{j}\equiv 0\pmod{p^{l+1}}.\label{c-40}
\end{align}
Combining \eqref{c-36}, \eqref{c-39} and \eqref{c-40}, we find that for $m>l+1$,
\begin{align}
\sum_{k}{'} \frac{1}{k^2}C^{(2,s,t)}_{n,p,m,l}(k)
\equiv \sum_{k}{'}\frac{1}{k^2}C^{(2,s,t)}_{n,p,m,l+1}(k)\pmod{p^{m+1}}.\label{c-41}
\end{align}
By a repeated use of \eqref{c-41}, we have
\begin{align}
\sum_{k}{'} \frac{1}{k^2} C^{(2,s,t)}_{n,p,m,1}(k)
\equiv \sum_{k}{'}\frac{1}{k^2}C^{(2,s,t)}_{n,p,m,m-1}(k)\pmod{p^{m+1}}.\label{c-42}
\end{align}

Note that
\begin{align}
&\sum_{k}{'}\frac{1}{k^2}C^{(2,s,t)}_{n,p,m,m-1}(k)\notag\\[5pt]
&=\sum_{N} \sum_{\substack{\lfloor k/p^{m-1}\rfloor=N\\[3pt] \{k/p^{m-1}\}<p^{m-1}/2}}{'}\frac{1}{k^2} C^{(2,s,t)}_{n,p,m,m-1}(k)
+\sum_{N} \sum_{\substack{\lfloor k/p^{m-1}\rfloor=N\\[3pt] \{k/p^{m-1}\}>p^{m-1}/2}}{'}\frac{1}{k^2} C^{(2,s,t)}_{n,p,m,m-1}(k).\label{c-43}
\end{align}
By \eqref{b-7} and \eqref{b-8}, we have
\begin{align}
\sum_{\substack{\lfloor k/p^{m-1}\rfloor=N\\[3pt] \{k/p^{m-1}\}<p^{m-1}/2}}{'}\frac{1}{k^2}
\equiv \sum_{\substack{\lfloor k/p^{m-1}\rfloor=N\\[3pt] \{k/p^{m-1}\}>p^{m-1}/2}}{'}\frac{1}{k^2}\equiv 0\pmod{p^{m-1}}.\label{c-44}
\end{align}
Applying \eqref{b-2} and \eqref{c-44} to the right-hand side of \eqref{c-43}, we arrive at
\begin{align}
&\sum_{k}{'}\frac{1}{k^2}C^{(2,s,t)}_{n,p,m,m-1}(k)
\equiv \sum_{k}{'}\frac{1}{k^2} C^{(2,s,t)}_{n,p,m,m}(k)\notag\\[5pt]
&\times\left(1-2np\sum_{j=1}^{\lfloor k/p^{m-1}\rfloor}{'}\frac{1}{j}
+snp\sum_{j=1}^{\lfloor k/p^{m-1}\rfloor}{'}\frac{1}{j}+tnp\sum_{j=1}^{\lfloor 2k/p^{m-1}\rfloor}{'}\frac{1}{j}\right)\pmod{p^{m+1}}.\label{c-45}
\end{align}
By \eqref{b-7} and \eqref{b-8}, we have
\begin{align}
&\sum_{k}{'}\frac{1}{k^2} C^{(2,s,t)}_{n,p,m,m}(k)\notag\\[5pt]
&=\sum_{N}{n-1\choose N}^2{n+N\choose N}^s{2N \choose n}^t
\sum_{\substack{\lfloor k/p^{m}\rfloor=N\\[3pt] \{k/p^{m}\}<p^{m}/2}}{'}\frac{1}{k^2}\notag\\[5pt]
&+\sum_{N}{n-1\choose N}^2{n+N\choose N}^s{2N+1 \choose n}^t
\sum_{\substack{\lfloor k/p^{m}\rfloor=N\\[3pt] \{k/p^{m}\}>p^{m}/2}}{'}\frac{1}{k^2}\notag\\[5pt]
&\equiv \frac{1}{3}p^mB_{p-3}\sum_{k}{n-1\choose k}^2{n+k\choose k}^s{2k \choose n}^t(12k+7)\notag\\[5pt]
&-\frac{1}{3}p^mB_{p-3}\sum_{k}{n-1\choose k}^2{n+k\choose k}^s{2k+1 \choose n}^t(12k+5)\pmod{p^{m+1}}.
\label{c-46}
\end{align}
By \eqref{b-14} and \eqref{b-15}, we have
\begin{align}
&{\sum_{k}}'\frac{1}{k^2} C^{(2,s,t)}_{n,p,m,m}(k)\sum_{j=1}^{\lfloor k/p^{m-1}\rfloor}{'}\frac{1}{j}\notag\\[5pt]
&=\sum_{N}{n-1\choose N}^2{n+N\choose N}^s{2N \choose n}^t
{\sum_{\substack{\lfloor k/p^{m}\rfloor=N\\[3pt] \{k/p^{m}\}<p^{m}/2}}}'\frac{1}{k^2}\sum_{j=1}^{\lfloor k/p^{m-1}\rfloor}{'}\frac{1}{j}\notag\\[5pt]
&+\sum_{N}{n-1\choose N}^2{n+N\choose N}^s{2N+1 \choose n}^t
{\sum_{\substack{\lfloor k/p^{m}\rfloor=N\\[3pt] \{k/p^{m}\}>p^{m}/2}}}'\frac{1}{k^2}\sum_{j=1}^{\lfloor k/p^{m-1}\rfloor}{'}\frac{1}{j}\notag\\[5pt]
&\equiv \frac{1}{3}p^{m-1} B_{p-3}\sum_{k}{n-1\choose k}^2{n+k\choose k}^s{2k \choose n}^t\notag\\[5pt]
&+\frac{1}{3}p^{m-1} B_{p-3}\sum_{k}{n-1\choose k}^2{n+k\choose k}^s{2k+1 \choose n}^t\pmod{p^m}.\label{c-47}
\end{align}
Similarly, by \eqref{b-16} and \eqref{b-17} we have
\begin{align}
&\sum_{k}{'}\frac{1}{k^2} C^{(2,s,t)}_{n,p,m,m}(k)\sum_{j=1}^{\lfloor 2k/p^{m-1}\rfloor}{'}\frac{1}{j}\notag\\[5pt]
&\equiv \frac{4}{3}p^{m-1} B_{p-3}\sum_{k}{n-1\choose k}^2{n+k\choose k}^s{2k \choose n}^t\notag\\[5pt]
&+\frac{4}{3}p^{m-1} B_{p-3}\sum_{k}{n-1\choose k}^2{n+k\choose k}^s{2k+1 \choose n}^t\pmod{p^m}.\label{c-48}
\end{align}
Combining \eqref{c-42} and \eqref{c-45}--\eqref{c-48} gives
\begin{align}
&{\sum_{k}}' \frac{1}{k^2} C^{(2,s,t)}_{n,p,m,1}(k)\notag\\
&\equiv \frac{1}{3}p^mB_{p-3}\sum_{k}{n-1\choose k}^2{n+k\choose k}^s{2k \choose n}^t(12k+7+sn-2n+4tn)\notag\\
&+\frac{1}{3}p^mB_{p-3}\sum_{k}{n-1\choose k}^2{n+k\choose k}^s{2k+1 \choose n}^t(-12k-5+sn-2n+4tn)\pmod{p^{m+1}}.\label{c-49}
\end{align}

It follows from \eqref{c-25}, \eqref{c-31} and \eqref{c-49} that
\begin{align*}
&Z^{(2,s,t)}_{m,n,p}\\
&\equiv \frac{1}{3}p^{3m}B_{p-3}\sum_{k}{n-1\choose k}^2{n+k\choose k}^s{2k \choose n}^tn^2(12k+7+sn-2n+4tn)\\
&+\frac{1}{3}p^{3m}B_{p-3}\sum_{k}{n-1\choose k}^2{n+k\choose k}^s{2k+1 \choose n}^tn^2(-12k-5+sn-2n+4tn)\\
&+ \frac{1}{2}p^{3m}B_{p-3}\sum_{k}{n-1\choose k}^2
{n+k\choose k}^s{2k+1 \choose n}^t n^3(t+2s+4) \\[5pt]
&-\frac{1}{2}p^{3m}B_{p-3}\sum_{k}{n-1\choose k}^2
{n+k\choose k}^s{2k \choose n}^t n^3(t+2s+4)\pmod{p^{3m+1}}.
\end{align*}

\vskip 5mm \noindent{\bf Acknowledgments.}
This work was supported by the National Natural Science Foundation of China (grant 12171370).

\end{document}